\documentclass{amsart}
\usepackage{amsmath, amssymb}
\usepackage[mathscr]{euscript}

\input xy
\xyoption{all}

\def\aut{\operatorname{Aut}}

\def\dcl{\operatorname{dcl}}

\def\tp{\operatorname{tp}}

\def\C{\mathcal{C}}

\def\id{\operatorname{id}}
\def\L{\mathcal{L}}
\def\acl{\operatorname{acl}}
\def\alg{\operatorname{alg}}
\def\acf{\operatorname{ACF}}
\def\acfa{\operatorname{ACFA}}

\def\dcf{\operatorname{DCF}}

\def\zloc{\operatorname{Z-loc}}
\def\trdeg{\operatorname{tr.deg.}}
\def\d{\operatorname d}
\def\dom{\operatorname{dom}}

\def\AA{\mathbb{A}}
\def\QQ{\mathbb{Q}}
\def\ZZ{\mathbb{Z}}
\def\KK{\mathbb{K}}

\def\Ind#1#2{#1\setbox0=\hbox{$#1x$}\kern\wd0\hbox to 0pt{\hss$#1\mid$\hss}
\lower.9\ht0\hbox to 0pt{\hss$#1\smile$\hss}\kern\wd0}
\def\ind{\mathop{\mathpalette\Ind{}}}
\def\Notind#1#2{#1\setbox0=\hbox{$#1x$}\kern\wd0\hbox to 0pt{\mathchardef
\nn=12854\hss$#1\nn$\kern1.4\wd0\hss}\hbox to
0pt{\hss$#1\mid$\hss}\lower.9\ht0 \hbox to
0pt{\hss$#1\smile$\hss}\kern\wd0}
\def\nind{\mathop{\mathpalette\Notind{}}}

\newtheorem*{theorem*}{Theorem}

\newtheorem{theorem}{Theorem}[section]
\newtheorem{proposition}[theorem]{Proposition}
\newtheorem{lemma}[theorem]{Lemma}
\newtheorem{corollary}[theorem]{Corollary}

\theoremstyle{definition}
\newtheorem{definition}[theorem]{Definition}

\theoremstyle{remark}

\begin{document}

\title[Model theory and differential-algebraic geometry]{Six lectures on model theory and differential-algebraic geometry}
\author{Rahim Moosa}
\address{Rahim Moosa\\
University of Waterloo\\
Department of Pure Mathematics\\
200 University Avenue West\\
Waterloo, Ontario \  N2L 3G1\\
Canada}
\email{rmoosa@uwaterloo.ca}

\date{\today}
\maketitle

\setcounter{tocdepth}{1}
\tableofcontents

%\bigskip
%\section{Introduction}

What follows is a write-up of some lectures I gave in the Fall of 2021 at the Fields Institute in Toronto, as part of the Thematic Programme on Trends in Pure and Applied Model Theory.
The actual lectures were given in four 90 minute instalments, but I present them here as six lectures because the material organised itself better that way.
I have taken other liberties as well: corrected errors, filled gaps, and improved exposition.

The goal of this module was to give a quick introduction to the model theory of differential fields that puts differential-algebraic geometry at the center.
As such, fundamental algebraic and model theoretic aspects of the subject, that would normally form the core of such a course, are entirely omitted.
Instead, I have tried to keep my focus on the birational geometry of algebraic vector fields, and more generally $D$-varieties in the sense of Buium~\cite{Buium}.
Applications of model theory to differential-algebraic geometry is an active area of research, and the approach I have taken here is meant to both whet the student's appetite and prepare them for work in the subject.

The style of the lectures are rather informal, lacking in both rigour and detail.
While I make every attempt to explain the central ideas, I leave many proofs to the student and I give essentially no references.
Moreover, I am content to take as black boxes many important theorems, especially when I have nothing to add to their exposition as it already exists in the literature.

There are a number of resources out there on the model theory of differential fields that the reader can consult for a more thorough, and more traditional, introduction.
Among them let me only mention Dave Marker's very influential treatment of the subject in~\cite{mtof}, and Anand Pillay's chapter in~\cite{lamt}.
(It is maybe interesting to note that the latter is based on another, much longer, course on differential fields held during another thematic programme at the Fields Institute some 25 years earlier.)
These articles also contain the references that I have generally omitted, coming, for example, from the differential algebra literature.
There is very little overlap between the lecture notes I am presenting here and either of these precedents.

I am grateful to the Fields Institute for hosting what was, for me personally, a very productive and stimulating programme.
It was especially successful as an early effort to return to normal in-person academic life and collaborative work.

\bigskip
\section{From geometry to algebra to model theory}

\noindent
My goal in this first lecture is to describe how the classical geometry-algebra correspondence gives rise to the model theory of differentially closed (and indeed difference closed) fields.

But first, as a review and to set the stage, let us consider the model theory of pure algebraically closed fields from this approach.
At its core, and speaking very loosely, the geometry-algebra correspondence is between a geometric space and the algebra of functions on the space.
In particular, to affine $n$-space $\AA^n$ over a field~$k$ of characteristic zero we associate the polynomial ring $k[x]$ in the $n$ variables $x=(x_1,\dots,x_n)$, with coefficients in $k$.
These we naturally view as functions on $\AA^n$ in the sense that we can evaluate a polynomial at any $n$-tuple from any field extension of $k$.
To an algebraic subvariety $V\subseteq\AA^n$ over $k$ we associate the co-ordinate ring $k[V]=k[x]/I(V)$, where $I(V)$ is the ideal of polynomials vanishing on $V$.
That is, we restrict the polynomials to functions on $V$, identifying polynomials if they agree on $V$ -- or rather on $V(K)$ for all field extensions $K\supseteq k$.
So the algebraic counterpart of the variety $V$ is the finitely generated $k$-algebra $k[V]$, which is an integral domain if $V$ is irreducible.
In fact, every finitely generated integral $k$-algebra, $R$, arises in this way.
Indeed, fixing generators $a_1,\dots,a_n\in R$ we obtain a surjective $k$-algebra homomorphism $\phi:k[x]\to R$ taking each $x_i$ to~$a_i$, so that $R$ is isomorphic to $k[x]/\ker\phi$, and letting $V\subseteq\AA^n$ be the subvariety determined by setting the polynomials in $\ker\phi$ to zero, we have exhibited $R$ as the co-ordinate ring of $V$.
So, if you are interested in the geometry of (embedded) irreducible affine algebraic varieties over $k$ then you should study finitely generated integral $k$-algebras.
A model-theorist will, of course, do so in the natural language $\L_k=\{0,1,+,-,\times,(\lambda_r)_{r\in k}\}$ of $k$-algebras.
Here $\lambda_r$ denotes the unary function symbol to be interpreted as sclar multiplication by $r\in k$.
Being finitely generated is not axiomatisable, but we can at least consider the (universal) theory $T$ of integral $k$-algebras in this language.
To better understand $T$ we should look to the {\em existentially closed models}: those $M\models T$ with the property that whenever a system of polynomial equations and inequations (i.e., a conjunction of atomic and negated atomic formulas) with coefficients from $M$ has a realisation in some model of $T$ extending $M$ then it already has a realisation in~$M$.
The class of existentially closed models of $T$ is itself axiomatisable; given by the theory $\acf_k$ of algebraically closed fields extending $k$.
It is thus that the study of affine algebraic varieties over~$k$ leads to the first order theory $\acf_k$.
Moreover, we recover the co-ordinate rings we were originally interested in as precisely the finitely generated substructures of the models of $\acf_k$.

The path we have just described, from geometry to algebra to model theory, and then back again, is the template for the expansion of algebraic geometry that we now study.
We will add additional structure to the algebraic varieties that we consider, right at the beginning of the above process, and then trace where this leads to algebraically and model-theoretically.
That additional structure is a vector field.
Namely, we are interested in algebraic varieties, $V\subseteq\AA^n$, equipped with a polynomial function that picks out, for each point $v\in V$, a tangent vector to~$V$ at~$v$ in~$\AA^n$.
Let's make this precise:

\begin{definition}
\label{def:avf}
Suppose $V\subseteq\AA^n$ is an affine variety over $k$.
By the {\em tangent bundle} to $V$ we mean the subvariety $TV\subseteq\AA^{2n}$ over $k$ defined, in co-ordinates $x=(x_1,\dots,x_n)$ and $y=(y_1,\dots,y_n)$, by
\begin{eqnarray*}
f(x)&=& 0\\
\sum_{i=1}^n\frac{\partial f}{\partial x_i}y_i&=&0
\end{eqnarray*}
for all $f\in I(V)$.
The projection onto the $x$ co-ordinates gives us a surjective morphism $\pi:TV\to V$.

An {\em algebraic vector field on $V$} is then a morphism $s:V\to TV$ over $k$ which is a section to $\pi$, that is, $\pi\circ s=\id_V$.
We will also sometimes refer to the pair $(V,s)$ as an {\em algebraic vector field} over $k$.
\end{definition}

For any field extension $K\supseteq k$, and any point $v\in V(K)$, the fibre of $\pi:TV\to V$ over~$v$ is a linear subspace of $\AA^n$ defined over the field $k(v)$.
Staring at the equations, we see that it is in fact the familiar {\em tangent space to $V$ at $v$} which we denote by~$T_vV$.
The algebraic vector field $s:V\to TV$ is given by $s=(\id_V, s_1,\dots, s_n)$ where $s_1,\dots s_n\in k[V]$ and we have that $(s_1(v),\dots,s_n(v))\in T_vV(K)$.

I hope it is clear that algebraic vector fields are inherently of interest, and I will make no effort to justify this claim.
In any case, we will take them as our basic geometric objects of study.
Our first question is: what is the algebraic counterpart to $(V,s)$?
In other words, what algebraic structure does $s$ induce on the co-ordinate ring $k[V]$?
The answer (you will have guessed) is a derivation.

\begin{definition}
\label{def:der}
A {\em derivation} on a (commutative) ring $R$ is a function $\delta:R\to R$ that is additive, namely satisfying $\delta(a+b)=\delta(a)+\delta(b)$, and satisfies the Leibniz rule $\delta(ab)=\delta(a)b+a\delta(b)$, for all $a,b\in R$.
By the {\em constants} of $(R,\delta)$ we mean $R^\delta:=\{a\in R:\delta(a)=0\}$.
\end{definition}

Note that the constants $R^\delta$ form a subring.
Assuming $R$ is nontrivial we must have that $R^\delta$ contains the integers: $\delta(1)=\delta(1\cdot 1)=\delta(1)1+1\delta(1)=2\delta(1)$ forces $\delta(1)=0$.
Note also that if $R$ is a $k$-algebra, then a derivation $\delta$ on $R$ will be $k$-linear if and only if $k\subseteq R^\delta$.
Indeed,  $k$-linearity forces $\delta(\lambda)=\delta(\lambda 1)=\lambda\delta(1)=0$ for all $\lambda\in k$, and conversely, $k\subseteq R^\delta$ implies $\delta(\lambda a)=\delta(\lambda)a+\lambda\delta(a)=\lambda\delta(a)$ for all $a\in R$.

\begin{proposition}
\label{prop:avf-der}
Suppose $V\subseteq\AA^n$ is a subvariety over $k$.
If  $s=(\id_V,s_1,\dots,s_n)$ is an algebraic vector field on $V$ over $k$ then there is a unique $k$-linear derivation $\delta_s$ on $k[V]$ such that $\delta_s(x_i+I(V))=s_i$ for all $i=1,\dots,n$.

Moreover, every $k$-linear derivation on $k[V]$ is of the form $\delta_s$ for some algebraic vector field $s$ on $V$ over $k$.
\end{proposition}

\begin{proof}
First we define $\delta_s$ on the polynomial ring $k[x]$.
Write each $s_i=g_i+I(V)$ for some $g_i\in k[x]$.
Then there is a unique $k$-linear derivation on $k[x]$ satisfying $\delta_s(x_i)=g_i$.
Indeed, uniqueness is clear because $k$-linearity and the Leibniz rule ensure that a $k$-linear derivation on a $k$-algebra is determined by its action on generators.
For existence, I give you the formula and leave it to you to check that it works:
\begin{equation}
\label{eqn:diffpoly}
\delta_sf:=\sum_{i=1}^n\frac{\partial f}{\partial x_i}g_i
\end{equation}
for each $f\in k[x]$.

The next thing to observe is that $I(V)$ is a {\em $\delta_s$-ideal}: it is closed under the action of $\delta_s$.
Here we use that $s:V\to TV$.
Indeed, for every field extension $K\supseteq k$, and every point $v\in V(K)$, we have that $s(v)\in TV(K)$, and hence
$$\sum_{i=1}^n\frac{\partial f}{\partial x_i}(v)s_i(v)=0$$
for all $f\in I(V)$.
Since the above identity holds for all $K$-points of $V$ as we range over all field extensions $K$, and since $s_i=g_i+I(V)$, we get that
$$\sum_{i=1}^n\frac{\partial f}{\partial x_i}g_i\in I(V).$$
By~(\ref{eqn:diffpoly}), $\delta_sf\in I(V)$ for all $f\in I(V)$, as desired.

So $\delta_s:k[x]\to k[x]$ induces a $k$-linear derivation on $k[x]/I(V)=k[V]$, which we also denote by $\delta_s$.
This has the desired property of $\delta_s(x_i+I(V))=s_i$ for all $i=1,\dots,n$, by construction.
Uniqueness follows as before.

Finally, for the ``moreover" clause, suppose we begin with a $k$-linear derivation $\delta$ on $k[V]$.
Write $\delta(x_i+I(V))=:s_i$ for all $i=1,\dots,n$.
Using $k$-linearity and the Leibniz rule one verifies that for any $f\in k[x]$, any $K\supseteq k$ a field extension, and any $v\in V(K)$,
$$\delta (f+I(V))(v)=\sum_{i=1}^n\frac{\partial f}{\partial x_i}(v)s_i(v).$$
Applying this to $f\in I(V)$ we get that $(s_1(v),\dots,s_n(v))\in T_vV(K)$.
That is, $s:=(\id_vs_1,\dots,s_n)$ is an algebraic vector field on $V$.
That $\delta=\delta_s$ is clear from construction.
\end{proof}

We see therefore, using the geometry-algebra correspondence, that to study irreducible algebraic vector fields we should consider the class of finitely generated integral $k$-algebras equipped with a $k$-linear derivation.
The natural language for this is $\L_{k,\delta}=\{0,1,+,-,\times,(\lambda_r)_{r\in k},\delta\}$.
We still cannot express the property of being finitely generated as a $k$-algebra, but we do have a universal theory, $T_{k,\delta}$, of integral $k$-algebras equipped with a $k$-linear derivation.
It turns out that the class of existentially closed models of $T_{k,\delta}$ is also elementary; its theory, that of {\em differentially closed fields} which contain $k$ and where the derivation vanishes on $k$, is denoted by $\dcf_k$.
Just as the study of algebraic varieties lead us to $\acf$, the study of algebraic vector fields motivates the model-theoretic consideration of $\dcf$.

But what about going back again?
That is, how do we recover the differential co-ordinate rings $(k[V],\delta_s)$ of algebraic vector fields that we were originally interested in?
This time, studying the finitely generated substructures of models of $\dcf_k$ won't do the trick as these co-ordinate rings are outright finitely generated as $k$-algebras, and not just as differential $k$-algebras.
Nevertheless, the differential rings $(k[V],\delta_s)$ can be detected from $\dcf_k$, as the {\em finite dimensional} substructures of models.
More on this later.

In these lectures I will be focusing on $\dcf$.
However, while we have this template set-up, let us consider one other variant breifly.
Instead of expanding $V\subseteq\AA^n$ by an algebraic vector field, consider instead an {\em algebraic dynamical system} on $V$, i.e., a dominant morphism $\phi:V\to V$ over $k$.
Recall that being dominant means $\phi(V(K))\subseteq V(K)$ is Zariski dense, for any (equivalently some) algebraically closed field extension $K\supseteq k$.
Again, I take for granted the intrinsic interest in algebraic dynamical systems $(V,\phi)$.
To understand the algebraic counterpart of this geometric object we need to think about what $\phi$ induces on the co-ordinate $k[V]$.
Writing $\phi=(f_1,\dots,f_n)$ with $f_1,\dots,f_n\in k[x]$, we obtain a $k$-linear endomorphism $\phi^*:k[V]\to k[V]$ given by $\phi^*(g)=g(f_1,\dots,f_n)$.
That is, viewing $g$ as a function on $V$, $\phi^*(g)$ is obtained by pre-composing with $\phi$.
In diagrams:
\xymatrix{
V\ar[rd]_{\phi^*(g)}\ar[r]^\phi &V\ar[d]^g\\
&\AA^1}.
Indeed, the geometry-algebra correspondence is a functor and this is how it acts on endomorphism.
The fact that $\phi$ is dominant implies that $\phi^*$ is injective: if $\phi^*(g)=0$ then $g$ must vanish on the image of $\phi$, which being Zariski dense forces $g=0$.
The algebraic counterpart to $(V,\phi)$ is the {\em difference} $k$-algebra $(k[V],\phi^*)$.
And again, all such arise in this way.
That is, every injective $k$-linear endomorphism of $k[V]$ is of the form $\phi^*$ for some algebraic dynamical system $\phi$ on $V$ over $k$.
To study irreducible affine algebraic dynamical systems over~$k$ is, therefore, to study finitely generated integral $k$-algebras equipped with an injective $k$-linear endomorphism.
The latter are studied model theoretically by the universal theory $T_{k,\sigma}$, in the language $\L_{k,\sigma}=\{0,1,+,-,\times,(\lambda)_{r\in k},\sigma\}$, of integral $k$-algebras equipped with an injective $k$-linear endomorphism.
The existentially closed models of $T_{k,\sigma}$ are axiomatisable, by the theory $\acfa_k$ of {\em difference closed fields} extending $(k,\id)$.
So, the same path leading from varieties to $\acf$, and from algebraic vector fields to $\dcf$, takes us from algebraic dynamical systems to $\acfa$.

\bigskip
\section{$D$-varieties}

\noindent
In the first lecture we worked exclusively over a fixed field $k$ on which our derivations were assumed to act trivially.
This is the so-called {\em autonomous} situation.
But the constraint is somewhat artificial, at least form the model-theoretic point of view.
In any case, even if we are interested primarily in the autonomous case, we will sometimes have to take base extensions to nontrivial differential fields in order to see the full geometric picture.
So, let us now fix a differential field $(k,\delta)$ of characteristic zero.

\begin{definition}
\label{def:dvar}
Suppose $V\subseteq\AA^n$ is an affine variety over $k$.
By the {\em prolongation} of $V$ over $(k,\delta)$ we mean the subvariety $\tau V\subseteq\AA^{2n}$ over $k$ defined, in co-ordinates $x=(x_1,\dots,x_n)$ and $y=(y_1,\dots,y_n)$, by
\begin{eqnarray*}
f(x)&=& 0\\
f^\delta(x)+\sum_{i=1}^n\frac{\partial f}{\partial x_i}y_i&=&0
\end{eqnarray*}
for all $f\in I(V)$.
Here $f^\delta$ denotes the polynomial obtained form $f$ by applying~$\delta$ to its coefficients.
The projection onto the $x$ co-ordinates gives us a surjective morphism $\pi:\tau V\to V$.

An affine {\em $D$-variety over $(k,\delta)$} is a pair $(V,s)$ where $V\subseteq\AA^n$ is an affine variety over $k$ and $s:V\to \tau V$ is a morphism over $k$ which is a section to $\pi$.
\end{definition}

Prolongations are the appropriate modification of tangent bundles in the presence of a derivation on the base field.
In particular, if $\delta=0$ on $k$ then $\tau V=TV$.
In general,  for each $v\in V$, the fibre of the prolongation, $\tau_vV$, is a coset in $\AA^n$ of the tangent space $T_vV$.
The point here is that when working in the possibly nonautonomous case, the basic geometric objects of interest are $D$-varieties rather than vector fields.

Where do the equations for the prolongation space come from?
The key differential-algebraic fact is the following elementary computation, the autonomous case of which was implicit in the proof of Proposition~\ref{prop:avf-der}.

\begin{lemma}
\label{lem:differentiatepoly}
Suppose $(R,\delta)$ is a differential ring, $f\in R[x]$ is a polynomial in $x=(x_1,\dots,x_n)$, and $a=(a_1,\dots,a_n)\in R^n$.
Then
$$\delta(f(a))=f^\delta(a)+\sum_{i=1}^n\frac{\partial f}{\partial x_i}(a)\delta(a_i).$$
\end{lemma}

\begin{proof}
We give a sketch, leaving the computations to the reader.
First prove the result, by induction on total degree, in the case when  $f$ is a monomial (and so, in particular $f^\delta=0$).
Next, consider a polynomial with only one term, say $f=bg$ where $g$ is a monomial and $b\in R$.
The result follows for such $f$ using the Leibniz rule and the case of monomials.
Finally, for the general case, note that all the operators involved in the desired identity -- namely $f\mapsto \delta (f(a))$,  $f\mapsto \frac{\partial f}{\partial x_j}(a)$ and $f\mapsto f^\delta(a)$ -- are additive in $f$.
So, as every polynomial is a sum of polynomials of the form already dealt with, the Lemma is proven.
\end{proof}

\begin{corollary}
\label{cor:nabla2prolong}
Suppose $(V,s)$ is a $D$-variety over $(k,\delta)$ and $(K,\delta)\supseteq (k,\delta)$ is a differential field extension.
Then $a\mapsto(a,\delta(a))$ defines a map $\nabla:V(K)\to\tau V(K)$.
\end{corollary}

\begin{proof}
It suffices to show that $f^\delta(a)+\sum_{i=1}^n\frac{\partial f}{\partial x_i}(a)\delta(a_i)=0$ for each $f\in I(V)\subseteq k[x]$ and $a\in V(K)$.
But, as $f(a)=0$, and hence $\delta(f(a))=0$, this is just Lemma~\ref{lem:differentiatepoly} applied to $(R,\delta)=(K,\delta)$.
\end{proof}

Now, the geometry-algebra correspondence yields a bijective correspondence between $D$-varieties over $(k,\delta)$ and finitely generated reduced $k$-algebras equipped with a derivation extending $\delta$.
That is, we have the following generalisation of Proposition~\ref{prop:avf-der}.

\begin{proposition}
\label{prop:dvar-der}
Suppose $V\subseteq\AA^n$ a subvariety over $k$.
If  $s=(\id_V,s_1,\dots,s_n)$ is a $D$-variety structure on $V$ over $k$ then there is a unique derivation $\delta_s$ on $k[V]$ extending $\delta$ on $k$ such that $\delta_s(x_i+I(V))=s_i$ for all $i=1,\dots,n$.

Moreover, every derivation on $k[V]$ extending $\delta$ on $k$ is of the form $\delta_s$ for some $D$-variety $s$ on $V$ over $k$.
\end{proposition}

\begin{proof}
The proof of Proposition~\ref{prop:avf-der} readily generalises to this nonautonomous setting, and I leave the details to you.
\end{proof}

Fix a $D$-variety $(V,s)$ over $(k,\delta)$.

\begin{definition}
\label{def:dsubvar}
A {\em $D$-subvariety} of $(V,s)$ is an algebraic subvariety $W\subseteq V$ over $k$ such that $s$ restricts to a $D$-variety structure on $W$.
That is, $s\upharpoonright_W:W\to\tau W$.
\end{definition}

It is worth thinking about what this means algebraically.
We have a differential structure $\delta_s$ on $k[V]$ induced by $s$, given by Proposition~\ref{prop:dvar-der}, and we have an ideal $I(W)\subseteq k[V]$ that defines the subvariety $W$.

\begin{lemma}
\label{lem:dsubvar-ideal}
A subvariety $W\subseteq V$ is a $D$-subvariety of $(V,s)$ if and only if $I(W)$ is a $\delta_s$-ideal of $k[V]$.
\end{lemma}

\begin{proof}
We are being a bit imprecise here.
If $V\subseteq\AA^n$ and $x=(x_1,\dots,x_n)$ are co-ordinates for affine $n$-space then $I(W)$ is an ideal of $k[x]$ containing $I(V)$.
But these are in bijective correspondence with ideals of $k[V]=k[x]/I(V)$, and it is in this sense that we view $I(W)$ in $k[V]$.

Write $s=(\id_V(s_1,\dots,s_n)$ and each $s_i=g_i+I(V)$ for some $g_i\in k[x]$.
Suppose $f\in I(W)$.
Then, by the constructuion of $\delta_s$ and Lemma~\ref{lem:differentiatepoly}, we have
$$\delta_s (f+I(V))=f^\delta+\sum_{i=1}^n\frac{\partial f}{\partial x_i}fg_i+I(V).$$
So $I(W)$ is a $\delta_s$-ideal if and only if the right-hand-side is in $I(W)$ for all $f\in I(W)$.
That is, we need to show that the following are equivalent:
\begin{itemize}
\item[(i)]
$W$ is a $D$-subvariety,
\item[(ii)]
for all $f\in I(W)$, $K\supseteq k$ a field extension, and all $a\in W(K)$,
$$f^\delta(a)+\sum_{i=1}^n\frac{\partial f}{\partial x_i}f(a)g_i(a)=0.$$
\end{itemize}
But, by definition $\tau W$ is defined by $f^\delta+\sum_{i=1}^n\frac{\partial f}{\partial x_i}fy_i$  as $f$ varies in $I(W)$, and the $g_i$ define $s$.
So (ii) is precisely expressing that $s$ maps $W$ to $\tau W$, which is what it means to be a $D$-subvariety.
\end{proof}

Let us a say a word about base extension.
Model-theoretically we are used to passing to a larger parameter set without comment.
But in algebraic geometry it is more standard (and prudent) to distinguish notationally between the given variety $V$ over $k$ and its base extension $V_K$ to a field extension $K\supseteq k$.
The co-ordinate ring of $V_K$ is the tensor product $k[V]\otimes_kK$, and the ideal of $V_K$ is the extension ideal $I(V)K[V_K]$.
So $V_K$ is just the variety $V$ viewed as being over $K$ rather than $k$.
Suppose, now, that $(K,\delta)\supseteq (k,\delta)$ is a differential field extension.
Then prolongations commute with base extension: $\tau(V_K)=(\tau V)_K$.
This follows from the fact (which we have not stated nor proved) that in the defining equations of the prolongation of $V$, given in Definition~\ref{def:dvar}, we could have restricted to any fixed set of generators for $I(V)$ instead of ranging over all $f\in I(V)$.
I leave this as an exercise for you.
It follows that $(V_K,s_K)$ is naturally a $D$-variety over $(K,\delta)$.
We will usually continue the model-theoretic habit of dropping these subscripts, and simply view $(V,s)$ as a $D$-subvariety over $(K,\delta)$ as well.
This allows us to use terminology like ``a $D$-subvariety of $(V,s)$ over $K$" when what is really meant is a $D$-subvariety of $(V_K,s_K)$.

\bigskip
\section{Finite dimensional types}

\noindent
Now let us pass to a differentially closed field.
These were discussed (motivated) in the first lecture, but I will deviate slightly from the conventions established there.
We will use the simpler language $\L_\delta:=\{0,1,+,-,\times,\delta\}$ of differential rings, rather than, say, the language $\L_{\QQ,\delta}$ of differential $\QQ$-algebras.
And instead of working with the theory $T_{\QQ,\delta}$ of integral $\QQ$-algebras equipped with a (necessarilly $\QQ$-linear) derivation, we consider the (still universal) $\L_\delta$-theory $T_\delta$ of differential integral domains of characteristic zero.
Note that if $(R,\delta)\models T_\delta$, and $R_{\QQ}$ denotes the localisation at $\ZZ\setminus\{0\}$, then $\delta$ extends uniquely to $R_{\QQ}$ and $(R_{\QQ},\delta)\models T_{\QQ,\delta}$.
So this change is really very harmless.
In particular, $T_\delta$ and $T_{\QQ,\delta}$ have the ``same" existentially closed models.
We denote the theory of these existentially closed models by  $\dcf_0$, and call it the the theory of {\em differentially closed fields of characteristic zero}.
Note that  I never proved that the class of existentially closed models of $T_\delta$ (or $T_{\QQ,\delta}$ for that matter) is elementary; it is, but I forego all discussion of the axiomatisation here.
Instead we will use existential closedness directly whenever we want to establish any properties of the models of $\dcf_0$.
For example, let me point out that differentially closed fields are algebraically closed:

\begin{lemma}
\label{lem:dcfisacf}
If $(K,\delta)\models\dcf_0$ then $K$ is an algebraically closed field.
\end{lemma}

\begin{proof}
Despite the terminology, we have not yet observed that $K$ is a field.
So let's do that first.
To see that every nonzero $a\in K$ is invertible, just apply existential closedness to the formula $\phi(x)$ given by $xa=1$.
Indeed, the quotient rule -- itself an immediate consequence of the Leibniz rule -- gives us a unique extension of $\delta$ to the localisation of $K$ at $a$, which thus yields a model of $T_\delta$ extending $(K,\delta)$ in which $\phi(x)$ has a realisation.
Hence $\phi(x)$ is already realised in $(K,\delta)$, as desired.

As similar approach gives that $K$ is algebraically closed.
Let $K(a)$ be a simple algebraic extension of $K$ and $f$ the minimal polynomial of $a$ over $K$.
The there is a unique extension of $\delta$ to $K(a)$.
Indeed, uniqueness -- which we don't actually need in this proof -- is by Lemma~\ref{lem:differentiatepoly} which dictates that $\delta(a)$ must equal $-\frac{f^\delta(a)}{f'(a)}$.
In any case, existence takes a bit more work.
First, view the derivation $\delta$ on $K$ as being $K(a)$-valued.
Then use the freeness of the polynomial ring to extend $\delta$ to a derivation $\delta:K[x]\to K(a)$ by setting $\delta(x):=-\frac{f^\delta(a)}{f'(a)}$.
Then observe, using Lemma~\ref{lem:differentiatepoly}, that $\delta(f)=0$.
We thus obtain an induced derivation $\delta:K[x]/(f)\to K(a)$.
The natural identification of $K[x]/(f)$ with $K(a)$ gives us our desired extension of $(K,\delta)$ to $(K(a),\delta)$.
The formula $f(x)=0$ has a realisation in $(K(a),\delta)\models T_\delta$, and hence by existential closedness, in $(K,\delta)$.
This forces $K(a)=K$, and we have shown that $K$ is algebraically closed.
\end{proof}

My goal for the rest of this lecture is to show how $D$-varieties capture precisely the finite dimensional fragment of $\dcf_0$.
Let us fix from now on a sufficiently saturated model $(\KK,\delta)\models\dcf_0$.
By convention all differential fields we consider are differential subfields of $\KK$ of cardinality strictly less than the level of saturation, unless explicitly stated otherwise.
(In fact, $\dcf_0$ is $\omega$-stable and hence admits saturated models of arbitrary large cardinality, so we can assume $\KK$ is saturated in its own cardinality.)
Similarly all parameters sets are assume to be small unless explicitly stated otherwise.

Fix a $D$-variety $(V,s)$ over a differential field $(k,\delta)$.

\begin{definition}
\label{def:dpoint}
By a {\em $D$-point of $(V,s)$} we mean a point $a\in V(\KK)$ such that the singleton $\{a\}$ is a $D$-subvariety of $(V,s)$ over $(\KK,\delta)$.
We denote the set of all $D$-points by $(V,s)^\sharp$.
\end{definition}

Note that we are implicitly taking a base extension as discussed earlier; the singleton $\{a\}$ is really a $D$-subvariety of the base extension of $(V,s)$ to $(\KK,\delta)$.
In order to preserve our convention that all parameters sets be small, we could instead take the base extension of $(V,s)$ to the differential field generated by $a$ over~$k$, namely $k\langle a\rangle:=k(a,\delta a,\delta^2 a,\dots)$.
In any case, we have the following characterisation of $(V,s)^\sharp$.

\begin{lemma}
$(V,s)^\sharp=\{a\in V(\KK):s(a)=\nabla(a)\}$.
\end{lemma}

\begin{proof}
The lemma follows immediately from the definitions once we observe that, for any $a\in V(\KK)$, the prolongation of $\{a\}$, viewed as a subvariety of $V$ over $F:=k\langle a\rangle$, is precisely $\{\nabla(a)\}$.
And this can be checked directly:
Suppose $V\subseteq\AA^n$ and $a=(a_1,\dots,a_n)$.
Then the ideal of $\{a\}$ is generated by the linear polynomials $x_1-a_1, x_2-a_2,\dots,x_n-a_n$ in $F[x_1,\dots,x_n]$.
So $\tau\{a\}\subseteq\AA^{2n}$ is defined in the variables $(x,y)$ by the equations $x_j-a_j=0$, for each $j=1,\dots,n$, along with the equation
$\displaystyle (x_j-a_j)^\delta+\sum_{i=1}^n\frac{\partial(x_j-a_j)}{\partial x_i}y_i=0$.
But this latter is just $y_j-\delta(a_j)=0$.
So the prolongation of $\{a\}$ is the singleton $\{(a_1,\dots,a_n,\delta a_1,\dots,\delta a_n)=\nabla(a)\}$.
\end{proof}

It follows that $(V,s)^\sharp$ is a $k$-definable set in $(\KK,\delta)$.
Indeed, writing the section as $s=(\id_V,s_1,\dots,s_n)$, and each $s_i=g_i+I(V)$ for some $g_1,\dots, g_n\in k[x]$, we have that $(V,s)^\sharp$ is defined by the formula
$\displaystyle (x\in V)\wedge\bigwedge_{i=1}^n \big(\delta x_i=g_i(x)\big)$.
Note that this formula is of a particularly simple form.
First of all, it is quantifier-freee.
But this is not surprising as in fact $\dcf_0$ admits quantifier elimination.
Moreover, it is a conjunction of {\em $\delta$-polynomial} equations; namely polynomial equations over $k$ in the variables $x,\delta x, \delta^2x,\dots$.
Such definable sets are called {\em Kolchin closed} in analogy with the Zariski closed sets of algebraic geometry.
In addition, $(V,s)^\sharp$ is of {\em order}~$1$ in the sense that only the first derivative of the variables actually appear.

Next, assuming that $V$ is irreducible, we associate to $(V,s)$ a certain complete type in $(\KK,\delta)$.
Let $p(x)$ be the collection of formulae over $k$ asserting that $x$ is a $D$-point of $(V,s)$ not contained in any proper subvariety of $V$ over $k$.

I claim that $p(x)$ is consistent.
That is, given a proper subvariety $W\subsetneq V$ over~$k$, there is a $D$-point of $(V,s)$ in $U:=V\setminus W$.
Our assumption that $V$ is irreducible ensures that $U$ is a dense Zariski open subset of $V$.
As you might expect, we need to use existential closedness of $(\KK,\delta)$ to show that there is a $D$-point in $U$, or indeed that there are any $D$-points of $(V,s)$ at all.
This is done as follows:
Using Proposition~\ref{prop:dvar-der}, let $\delta_s$ be the derivation on $k[V]$ extending $\delta$ on $k$ that is induced by $s$.
Let $a:=x+I(V)$ be the generators of $k[V]=k[x]/I(V)$ coming form the variables of the ambient polynomial ring.
By definition of $\delta_s$ we have that $\delta_s(a_i)=s_i$ for all $i=1,\dots,n$, in the differential ring $(k[V],\delta_s)$.
That is, $a$ is a realisation in $(k[V],\delta_s)$ of the formula defining $(V,s)^\sharp$.
Moroever, $a$ also realises $x\notin W$ as $I(W)\supsetneq I(V)$ since $W$ is a proper subvariety of $V$.
That is, the formula defining $(V,s)^\sharp\cap U$ has a realisation in $(k[V],\delta_s)$.
Now,  as $V$ is irreducible, $k[V]$ is an integral domain, and hence $(k[V],\delta_s)$ is a model of $T_\delta$ extending $(k,\delta)$.
It follows by existential closedness that the formula defining $(V,s)^\sharp\cap U$ must have a realisation in $(\KK,\delta)$, as desired.

Next, I claim that $p(x)$ determines a complete type.
Using the fact that $\dcf_0$ admits quantifier elimination (whose proof I also forego!), it suffices to prove that if $a\models p(x)$ then $\tp_\L(a,\delta a,\delta^2a,\dots/k)$ is determined, where recall that $\L$ is just the language of rings (as opposed to the language $\L_\delta$ of differential rings).
But as $a$ is a $D$-point of $(V,s)$, the derivative $\delta a$ is given by polynomials in $a$ over $k$, and hence so is $\delta^\ell a$ for all $\ell\geq 1$ by Lemma~\ref{lem:differentiatepoly}.
So $\tp_\L(a,\delta a,\delta^2a,\dots/k)$ is determined by $\tp_\L(a/k)$.
By quantifer elimination in $\acf_0$, this is in turn determined by $\zloc(a/k)$, the {\em Zariski locus} of $a$ over $k$; namely, the smallest Zariski closed subset of $V$ over $k$ that contains $a$.
But as $p(x)$ ensures that $a$ is not contained in any proper sunbvariety of $V$ over $k$, we must have that $\zloc(a/k)=V$.
So $\tp_\L(a/k)$ is determined.

\begin{definition}
\label{def:dvar-gen}
We call this $p(x)$ the {\em generic type of $(V,s)$ over $k$}.
It is the type asserting that $x$ is a $D$-point of $(V,s)$ and that the Zariski locus of $x$ over $k$ is $V$.
We call a realisation of $p(x)$ a {\em generic $D$-point of $(V,s)$ over $k$}.
\end{definition}

We have thus associated a complete type to every irreducible $D$-variety.
Not all complete types arise in this way.
To see this, let us introduce dimension for types as follows:

\begin{definition}
\label{def:dimension}
Suppose $(k,\delta)$ is a differential field and $p\in S(k)$ is a complete type over $k$.
By the {\em dimension of $p$}, denoted by $\dim(p)$, we mean the non-decreasing sequence of non-negative integers
$$(\trdeg(a/k),\trdeg(a,\delta a/k),\trdeg(a,\delta a,\delta^2a/k),\dots)$$
where $a\models p$.
We also write this as $\dim(a/k)$.
If this sequence eventually stabilises then we say that $p$ is {\em finite dimensional} and we write $\dim(p)=d$ where $d$ is that eventual value of the sequence.
\end{definition}

Note that $p$ is finite dimensional if and only if for some (equivalently any) $a\models p$ we have that the differential field $k\langle a\rangle= k(a,\delta a,\delta^2 a,\dots)$ generated by $a$ over $k$ has finite transcendence degree, and in that case $\dim(p)=\trdeg(k\langle a\rangle/k)$.

\begin{lemma}
\label{lem:dimgentype}
If $(V,s)$ is an irreducible $D$-variety over $(k,\delta)$ then its generic type  over $k$ is finite dimensional and of dimension $\dim V$.
\end{lemma}

\begin{proof}
If $a$ is a generic $D$-point of $(V,s)$ over $k$ then $\trdeg(a/k)=\dim V$ since $V=\zloc(a/k)$.
On the other hand, $k(a,\delta a,\dots,\delta^\ell a)=k(a)$ for all $\ell\geq 1$ as $a$ is a $D$-point and hence $\delta a $ is a polynomial in $a$ over $k$.
So $\dim(a/k)=(\dim V,\dim V,\dots)$.
\end{proof}

In fact, every finite dimensional type arises in this way.
Well, at least up to interdefinability.
Here, we say that $p,q\in S(k)$ are {\em interdefinable} if for all (equivalently some) $a\models p$ there is $b\models q$ such that $\dcl(ka)=\dcl(kb)$.
Using quantifier-elimination for $\dcf_0$ one can show that $\dcl(ka)=k\langle a\rangle$, and hence finite dimensionality is an interdefinability invariant of complete types.
But note that dimension itself is {\em not} an interdefinability invariant: if $a\in\KK$ is differentially-transcendental over $k$ in the sense that $(a,\delta a,\delta^2a,\dots)$ is an algebraically independent sequence then $\dim(a/k)=(1,2,3,\dots)$ while $\dim (a,\delta a/k)=(2,3,4,\dots)$, though $\tp(a/k)$ and $\tp(a,\delta a/k)$ are interdefinable.

\begin{theorem}
\label{thm:fdt}
Suppose $(k,\delta)$ is a differential field.
Every finite dimensional complete type over $k$ is interdefinable with the generic type of an irreducible $D$-variety over $(k,\delta)$.
\end{theorem}

\begin{proof}
I will give only a sketch, in the special case of $1$-types, from which you will see how to proceed in general and in detail.
We are given a type $p=\tp(a/k)$ where $a\in\KK$ and $k\langle a\rangle$ is of finite transcendence degree over $k$.
It follows that, for some $\ell\geq 0$, $\delta^\ell a\in k(a,\delta a,\dots,\delta^{\ell-1}a)^{\alg}$. 
Let $P$ be the minimal polynomial of $\delta^\ell a$ over $k(a,\delta a,\dots,\delta^{\ell-1}a)$.
Since $P(\delta^\ell a)=0$, differentiating both sides we get from Lemma~\ref{lem:differentiatepoly} that
$$0=\delta(P(\delta^\ell a))=P^\delta(\delta^\ell a)+P'(\delta^\ell a)\delta^{\ell+1}a.$$
By minimality, $P'(\delta^\ell a)\neq 0$, and hence we get that
$$\delta^{\ell+1}a=-\frac{P^\delta(\delta^\ell a)}{P'(\delta^\ell a)}=f(a,\delta a,\dots,\delta^\ell a)$$
for some rational function $f\in k(x^{(0)},x^{(1)},\dots,x^{(\ell)})$.
Now let $b:=(a,\delta a,\dots,\delta^\ell a)$.
Then
$$\nabla (b)=(a,\delta a,\dots,\delta^\ell a, \delta a,\delta^2 a,\dots,\delta^{\ell+1} a)=(b,\pi_1(b),\pi_2(b),\dots,\pi_\ell(b), f(b))$$
where the $\pi_i$ are the $i$th co-ordinate projections.
That is, $\nabla(b)=s(b)$ where $s=(\id,\pi_1,\dots,\pi_\ell,f)$.
If we set $V:=\zloc(b/k)\setminus\operatorname{pole}(f)\subseteq\AA^{\ell+1}$, then $s:V\to \tau V$ is a regular section, $b$ is a generic $D$-point of $(V,s)$ over $k$, and $\tp(a/k)$ and $\tp(b/k)$ are interdefinable.
I am cheating, of course, because $V$ is not a closed subvariety of $\AA^{\ell+1}$, but by working with one more variable this can be remedied.
\end{proof}

The study of the birational geometry of $D$-varieties thus coincides with the model theory of the finite dimensional fragment of $\dcf_0$.

\bigskip
\section{Stability and independence}

\noindent
I have omitted proofs of some of the most fundamental properties of $\dcf_0$; in particular, that it exists (i.e., the fact that the existentially closed models of $T_\delta$ form an elementary class) and that it admits quantifier elimination.
This is largely because you can find the proofs elsewhere, and I had nothing to add.
In this lecture I want to discuss some further model theoretic properties, around stability, and this time I will give at least some proofs.

Before talking about stability, let me say a few words about the elimination of imaginaries, another important property that $\dcf_0$ enjoys.
This is the statement that every definable set $D$ has a {\em code}; that is, a finite tuple $e$ such that for all $\sigma\in \aut(\KK,\delta)$, $\sigma(D)=D$ if and only if $\sigma(e)=e$.
Equivalently, for some formula $\phi(x,y)$, $D=\phi(\KK,e)$ but $D\neq\phi(\KK,e')$ for any $e'\neq e$.
That is, a code for a definable set is a kind of minimal and canonical parameter.
Using quantifier elimination and the noetherianity of the Kolchin topology (both facts not proved here), it is not hard to reduce the verification that all definable sets have codes to showing that all Kolchin closed sets have codes.
Now, a Kolchin closed set is of the form
$$D=\{a\in\KK^n:\nabla_\ell(a):=(a,\delta a,\delta^2a,\dots,\delta^\ell a)\in V(\KK)\}$$
for some $\ell\geq 0$ and some algebraic subvariety $V\subseteq\AA^{(\ell+1)n}$.
Moreover, replacing $V$ by the Zariski closure of $\nabla_\ell(D)$, we may assume that $\nabla_\ell(D)$ is Zariski dense in $V$.
From algebraic geometry we know that $V$ has a minimal field of definition, say $L$, that $L$ is finitely generated, say $L=\QQ(e)$, and that $e$ is a code for $V(\KK)$ in the pure field structure $\KK$.
I claim that $e$ will in fact be a code for $D$ in $(\KK,\delta)$.
It is clear that $D$ is defined over $e$, and so it suffices to show that if $\sigma\in \aut(\KK,\delta)$ preserves $D$ then $\sigma(e)=e$.
But as $\sigma$ commutes with $\delta$, if $\sigma(D)=D$ then $\sigma(\nabla_\ell(D))=\nabla_\ell(D)$, and hence $\sigma$ preserves the Zariski closure $V$, which in turn forces $\sigma(e)=e$.

Now let us pass to stability.

\begin{theorem}
\label{thm:stability}
$\dcf_0$ is {\em $\omega$-stable}: there are only countably many types over countably many parameters.
\end{theorem}

\begin{proof}
It suffices to count $1$-types.
Actually, the proof sketch I gave of Theorem~\ref{thm:fdt} already suggests how to count the $1$-types over a differential field $(k,\delta)$, and we follow that suggestion now.
First of all, there is a unique $1$-type of a {\em differentially transcendental} element over $k$; that is, of an element $a\in\KK$ such that $(a,\delta a, \delta^2 a,\dots)$ is an algebraically independent sequence over $k$.
So it remains to count {\em differentially algebraic} $1$-types (namely those that are not differentially transcendental).

Suppose therefore that $a\in \KK$ is differentially algebraic over $k$.
Let $\ell\geq 0$ be least such that $\delta^\ell a\in k(a,\delta a,\dots,\delta^{\ell-1}a)^{\alg}$ and let $P(t)$ be the minimal polynomial of $\delta^\ell a$ over $k(a,\delta a,\dots,\delta^{\ell-1}a)$.
After clearing denominators we can write $P=g(a,\delta a,\dots,\delta^{\ell-1} a, t)$ where $g\in k[x^{(0)},x^{(1)},\dots,x^{(\ell)},t]$.
I claim that $\tp(a/k)$ is determined by the pair $(\ell, g)$.
But before proving this let us observe that $k\langle a\rangle=k(a,\delta a,\dots,\delta^{\ell}a)$.
Indeed, we saw in the proof of Theorem~\ref{thm:fdt} that $\delta^\ell a$ being algebraic over $k(a,\delta a,\dots,\delta^{\ell-1}a)$ implies that $\delta^{\ell+1} a$ is contained in $k(a,\delta a,\dots,\delta^{\ell}a)$.
In particular,  $\delta^{\ell+1} a$ is algebraic over $k(a,\delta a,\dots,\delta^{\ell}a)$ and hence
$$\delta^{\ell+2} a\in k(a,\delta a,\dots,\delta^{\ell+1}a)=k(a,\delta a,\dots,\delta^{\ell}a).$$
Iterating gives us the desired fact that $k\langle a\rangle=k(a,\delta a,\dots,\delta^{\ell}a)$.

Suppose now that $b\in \KK$ is differentially algebraic over $k$ and gives rise to the same data $(\ell, g)$.
I want to show that $\tp(a/k)=\tp(b/k)$.
To do so, I will exhibit a differential-field-isomorphism from $k\langle a\rangle$ to $k\langle b\rangle$ over $k$, that takes $a$ to $b$, and this will suffice by quantifier elimination.
But we know that
$$k\langle a\rangle=k(a,\delta a,\dots,\delta^{\ell}a)$$
and
$$k\langle b\rangle=k(b,\delta a,\dots,\delta^{\ell}b).$$
So it suffices to exhibit a field-isomorphism $\alpha: k(a,\delta a,\dots,\delta^{\ell}a)\to k(b,\delta b,\dots,\delta^{\ell}b)$, over $k$, which satisfies $\alpha(\delta^ia)=\delta^ib$ for all $i=0,\dots,\ell+1$.
(It is not a typo here that we have to check all the way up to $i=\ell+1$.)
Indeed, I am using the fact, which I leave to you to check, that if you have a field isomorphism between differential fields which commutes with the derivation on the field-generators of the domain, then it must be a differential-field-isomorphism.
(Hint: This too rests on the infinitely useful Lemma~\ref{lem:differentiatepoly}.)

First of all, $\tp_{\L}(a,\delta a,\dots,\delta^{\ell-1}a/k)=\tp_{\L}(b,\delta b,\dots,\delta^{\ell-1}b/k)$ as it is the field-type of an algebraically independent  $\ell$-tuple.
So we have a field isomorphism $\alpha:k(a,\delta a,\dots,\delta^{\ell-1}a)\to k(b,\delta b,\dots,\delta^{\ell-1}b)$, over $k$, taking $\delta^i a$ to $\delta^i b$ for all $i=0,\dots,\ell-1$.
We want to extend $\alpha$ to $\delta^\ell a$.
To do so, note that, by construction, $P=g(a,\delta a,\dots,\delta^{\ell-1}a,t)$ is an irreducible polynomial over $k(a,\delta a,\dots,\delta^{\ell-1}a)$ of which $\delta^\ell a$ is a root.
On the other hand, $\delta^\ell b$ is a root of $Q:=g(b,\delta b,\dots,\delta^{\ell-1}b,t)$, and $Q$ is the transform of $P$ by $\alpha$.
So we can extend $\alpha$ to
$$\alpha: k(a,\delta a,\dots,\delta^{\ell}a)\to k(b,\delta b,\dots,\delta^{\ell}b)$$
by sending $\delta^\ell a$ to~$\delta^\ell b$.

It remains only to check that $\alpha(\delta^{\ell+1}a)=\delta^{\ell+1}b$.
But this also follows from the proof of Theorem~\ref{thm:fdt}, where we saw that
$\delta^{\ell+1}a=-\frac{P^\delta(\delta^\ell a)}{P'(\delta^\ell a)}$ and $\delta^{\ell+1}b=-\frac{Q^\delta(\delta^\ell b)}{Q'(\delta^\ell b)}$.
I leave it to you to check that $\alpha$ takes $P^\delta$ to $Q^\delta$ and $P'$ to $Q'$.
Hence $\alpha$ takes $\delta^{\ell+1} a$ to $\delta^{\ell+1}b$, as desired.

We have proved that the differentially algebraic $1$-type $\tp(a/k)$ is determined by the pair $(\ell,g)$.
If $k$ is countable then there are only countably many possible such pairs, and hence only countably many differentially algebraic $1$-types over~$k$ (and only one differentially transcendental $1$-type).
As every countable set of parameters is contained in a countable differential field (namely the differential field it generates), it follows that over countably many parameters we have only countably many complete $1$-types.
That is, $\dcf_0$ is $\omega$-stable.
\end{proof}

What I like about the above proof is that it uses very little differential algebra.
In particular, somewhat unexpectedly, no use is made of the Ritt-Raudenbush basis theorem which says that every radical differential ideal in a differential polynomial ring over a differential field is finitely generated (as a radical differential ideal).
Differential algebra is a useful and beautiful subject, but it is interesting to note how little of it one really needs to do model theory in $\dcf_0$.

Once we have $\omega$-stability, the full machinery of geometric stability theory becomes available.
In particular we have the good behaviour of Shelah's nonforking independence.
I will not give the abstract definition of nonforking, but rather specialise to what it means in $\dcf_0$.

\begin{definition}
Given a tuple $a$ and subsets $B\subseteq A$ of $\KK$, we say that {\em $a$ is independent from $A$ over~$B$}, denoted by $\displaystyle a\ind_BA$, to mean that $\dim(a/\dcl(A))=\dim(a/\dcl(B))$.
In this case we also say that $\tp(a/A)$ a {\em does not fork over $B$} or that $\tp(a/A)$ is a {\em nonforking extension} of $\tp(a/B)$.
\end{definition}

Recall that the dimension of a type over a differential field was defined in~\ref{def:dimension} as a certain infinite sequence of nondecreasing integers, and it is as such that equality is to be understood above.
Namely, if $F:=\dcl(A)$ and $k:=\dcl(B)$, then $\displaystyle a\ind_BA$ if and only if $\trdeg(a,\delta a,\dots,\delta^\ell a/F)=\trdeg(a,\delta a,\dots,\delta^\ell a/k)$, for all $\ell\geq 0$.
It is also worth pointing out that this is equivalent to saying that $k\langle a\rangle$ is {\em algebraically disjoint}, in the sense of pure field theory, from $F$ over $k$.

Maybe the best way to see that this agrees with Shelah's nonforking is to prove that it satisfies the usual desired properties (nontriviality, finite character, automorphism invariance, symmetry, transitivity, existence of nonforking extensions, uniqueness of nonforking extensions over algebraically closed sets), using the analogous properties for algebraic independence in pure fields, and then use the fact that Shelah's nonforking in stable theories is characterised among all abstract independence relations by these properties.
But we do not go into that here.

%This will be stated in terms of {\em algebraic disjointness}: a subfield $K$ of $\KK$ is algebraically disjiont from another subfield $L$ of $\KK$ over a common subfield $F\subseteq K\cap L$ if there is a transcendence basis for $K$ over $F$ that remains transcendental over~$L$. In order for the following definition of independence to be meaningful, you should also know that, as a consequence of quantifier elimination, the model-theoretic algebraic closure of a set $A$, denoted $\acl(A)$, agrees with the field-theoretic algebraic closure of the differential field generated by $A$.

How does this specialise to generic types of $D$-varieties (and hence, by Theorem~\ref{thm:fdt}, to finite dimensional types)?
Here is the simple answer:

\begin{proposition}
\label{prop:indgentype}
Suppose $(V,s)$ is a $D$-variety over $(k,\delta)$ with $V$ absolutely irreducible\footnote{A variety is {\em absolutely irreducible} if its base extension to any  field extension is irreducible.}.
Let $a$ be a generic $D$-point of $(V,s)$ over $k$.
Then, for any differential field extension $F\supseteq k$, $a\ind_kF$ if and only if $a$ is a generic $D$-point of $(V,s)$ over~$F$.
\end{proposition}

\begin{proof}
By Lemma~\ref{lem:dimgentype}, the dimension of the generic type of $(V,s)$ over $k$ is $\dim V$, and so is the dimension of the generic type of $(V,s)$ over $F$.
(Note that we are using absolute irreducibility here so that we can apply Lemma~\ref{lem:dimgentype} to the base extension of $(V,s)$ to $F$, which remains irreducible.)
This proves the right-to-left direction.
For the converse, assume $a\ind_kF$.
As we already know that $a$ is a $D$-point of $(V,s)$, it remains to verify that the Zariski locus of $a$ over~$F$, say $W$, is equal to $V$.
But $W$ is a Zariski closed subset of $V$ over $F$ and it is of dimension $\trdeg(a/F)=\trdeg(a/k)=\dim V$.
Hence, by the irreducibility of $V$ over~$F$,  we must have $W=V$.
\end{proof}

The assumption in the above proposition that $V$ be absolutely irreducible may seem a little unnatural.
At the very least we should give a model-theoretic explanation of absolute irreducibility:

\begin{proposition}
\label{prop:statgentype}
Suppose $(V,s)$ is an irreducible $D$-variety over $(k,\delta)$.
Then $V$ is absolutely irreducible if and only if the generic type of $(V,s)$ over $k$ is {\em stationary}: it has a unique nonforking extension to every set containing $k$.
\end{proposition}

\begin{proof}
If $V$ is absolutely irreducible, then by Proposition~\ref{prop:indgentype} the nonforking extension of its generic type, $p$, over $k$, to a set $A\supseteq k$, is precisely the generic type of $(V,s)$ over $F:=\dcl(A)$.
Conversely, suppose the generic type $p$ of $(V,s)$ over $k$ is stationary and let $F:=k^{\alg}$.
For absolute irreducibility it suffices to verify that $V$ is irreducible over $F$.
Note that $F$ is also a differential subfield of $(\KK,\delta)$, since, as observed in the proof of Lemma~\ref{lem:dcfisacf}, if $a\in k^{\alg}$ with minimal polynomial $f$ then $\delta(a)=-\frac{f^\delta(a)}{f'(a)}$.
Now, each irreducible component $W$ of $V$ over $F$ is of dimension $\dim V$.
I claim that it suffices to check that $W$ is a $D$-subvariety of $V$ over $F$.
Indeed, in that case the generic type of $(W,s\upharpoonright_W)$ over $F$ would be the nonforking extension of $p$ to $F$, by Proposition~\ref{prop:indgentype}, and hence all irreducible components of $V$ over $F$ would share a generic type, forcing $V$ itself to be irreducible over $F$.

So we have reduced to showing that the irreducible components of $V$ over $F$ are $D$-subvarieties.
This is of independent interest.
Recall, from algebraic geometry, that  the irreducible components of a variety correspond to the minimal prime ideals containing the ideal of that variety.
By Lemma~\ref{lem:dsubvar-ideal}, we need to show that those minimal prime ideals are $\delta_s$-ideals.
This follows from the following fundamental (but easy) fact of differential algebra: {\em If $(F,\delta)$ is a differential field and $R$ is a finitely generated $F$-algebra equipped with a $\delta$-ring structure extending $(F,\delta)$, and $I$ is a radical $\delta$-ideal of $R$, then every minimal prime ideal containing $I$ is also a $\delta$-ideal.}
To see this, let $I=P_1\cap\cdots\cap P_\ell$ be the prime decomposition of $I$.
We need to observe that each of the $P_j$ is a $\delta$-ideal.
That is, fixing $j=1,\dots,\ell$, and $a\in P_j$, we need to show that $\delta(a)\in P_j$.
For each $r\neq j$ choose $b_r\in P_r\setminus P_j$ and let $b$ be the product of all the $b_r$'s.
Then $ab\in I$, and hence $\delta(a)b+a\delta(b)=\delta(ab)\in I\subseteq P_j$.
Since $a\in P_j$ this implies $\delta(a)b\in P_j$.
But by construction $b\notin P_j$, and hence $\delta(a)\in P_j$ as desired.
\end{proof}

\bigskip
\section{Around the constants}

\noindent
By the {\em field of constants} we mean $\C:=\{a\in\KK:\delta(a)=0\}$.
Note that $\C$ is the set of $D$-points of the trivial vector field on the affine line; that is, $\C=(\AA^1,0)^\sharp$.
Here $0$ here denotes the zero vector field on the affine line given by $a\mapsto (a,0)$.
It is also easy to verify that $\C$ is a subfield of $\KK$.
Let us observe that it is algebraically closed:
If $f(t)$ is the minimal polynomial over $\C$ of an element $a\in\KK$ that is algebraic over $\C$, then by differentiating $0=f(a)$ we get $0=f'(a)\delta(a)$ as $f^\delta=0$.
But $f'(a)\neq 0$ by minimality, so that we must have $\delta(a)=0$.
This shows that $\C$ is relatively algebraically closed in $\KK$.
But we have seen that $\KK$ is algebraically closed (Lemma~\ref{lem:dcfisacf}), and hence so is $\C$.

\begin{proposition}
\label{prop:pureconst}
The constants form a {\em stably embedded pure} algebraically closed field.
That is, if $D\subseteq\C^n$ is definable in $(\KK,\Delta)$ then it is definable in $(\C,+,\times)$.
\end{proposition}

\begin{proof}
Note that for us ``definable" means ``definable with parameters".
So the first step is to show that $D$ is definable in $(\KK,\Delta)$ with parameters from $\C$.
We will use the fact (not proven in this course) that all types of a stable theory are definable.
Suppose $D$ is defined by $\phi(a,y)$ where $\phi(x,y)$ is a $\L_\delta$-formula and $a\in\KK^m$, and let ${\bf p}(x):=\tp(a/\C)$.
I use boldface here because $\bf p$ is not a proper type according to our conventions where we are only to allow parameter sets that are of cardinality less then the level of saturation.
But it is OK, such {\em global types} are also definable.
In particular, we have a formula $\d_{\bf p}\phi(y)$ over $\C$ such that
$\phi(x,b)\in\bf p$ if and only if $\models\d_{\bf p}\phi(b)$
for all $b\in\C^n$.
But this implies that, for all $b\in\KK^n$,
$$b\in D\ \iff\ \models \d_{\bf p}\phi(b)\wedge(\delta b=0)$$
which shows that $D$ is definable with parameters in $\C$.

Next, easy manipulations using quantifier elimination reduces us further to the case when $D$ is defined by a conjunction $(\delta y=0)\wedge (P(y,\delta y,\dots,\delta^\ell y)=0)$ where $P$ is a polynomial with coefficients in $\C$.
But that formula is equivalent to
$$(\delta y=0)\wedge (P(y,0,\dots,0)=0),$$
and the latter defines a set which is clearly definable in $(\C,+,\times)$.
\end{proof}

In this way, pure algebraic geometry lives definably in the finite dimensional fragment of $\dcf_0$, precisely as the induced structure on the constants.
(Notice that this is saying something different, more meaningful, than that $(\KK,\Delta)$ is an expansion of the pure algebraically closed field $\KK$.)
The model theorist's approach can now be described as follows:
Understand the fine structure of a finite dimensional type by studying its relationship to the field of constants.
In other words, study differential-algebraic geometry in relation to the algebraic geometry living therein.

Such an approach will, of necessity, say nothing about those types that have no definable relationship to the constants.
Here we have to be careful about what we might mean by ``having no definable relation to the constants", in particular with respect to parameters.
Consider, for example, the set $D\subseteq\KK$ defined by the equation $\delta(x)=1$.
Note that $D=(\AA^1,1)^\sharp$ where $1$ here denotes the constant vector field on the affine line given by $a\mapsto (a,1)$.
The elements of $D$ are independent of the constants over the empty set.
Indeed, if $d\in D$ and $F\subseteq\C$ is any subfield then, as $d\notin F^{\alg}\subseteq\C$, $\dim(d/F)=1=\dim(d/\QQ)$.
It follows that $d\ind F$, which expresses the fact that there are no nontrivial $0$-definable relations between $D$ and $\C$.
On the other hand, if we fix a solution $d_0\in D$, and let $k:=\QQ\langle d_0\rangle$, then $D$ is definably isomorphic to~$\C$ over~$k$ as $D=d_0+\C$.
In fact, we have an isomorphism of $D$-varieties $(\AA^1,1)\to(\AA^1,0)$ over $k$, given by translation by $-d_0$.
So parameters matter very much, leading to the following two natural implementations of ``having no definable relation to the constants".

\begin{definition}
\label{def:orthogonal}
Suppose $(k,\delta)$ is a differential field and $p\in S(k)$ is a complete stationary type.
We say that $p$ is {\em weakly orthogonal to $\C$}, denoted by $\displaystyle p\perp^w\C$, to mean that whenever $a\models p$ and $c$ is a tuple from $\C$ then $a\ind_kc$.
We say that $p$ is {\em orthogonal to $\C$}, denoted $p\perp\C$, if every nonforking extension of $p$ is weakly orthogonal to~$\C$.
\end{definition}

That is, $p\perp\C$ means that for any $B\supseteq k$, any $a\models p$ with $a\ind_kB$, and any tuple $c$ from $\C$, $a\ind_Bc$.
For example, what follows from the above discussion is that the type over $\QQ$ of any solution to $\delta(x)=1$ is weakly $\C$-orthogonal, but not $\C$-orthogonal.
In fact, $\delta(x)=1$ fails orthogonality to the constants in a particularly strong way:

\begin{definition}
Suppose $(k,\delta)$ is a differential field and $p\in S(k)$ is a complete stationary type.
We say that $p$ is ({\em almost}) {\em $\C$-internal} if there is $B\supseteq k$, $a\models p$ with $a\ind_kB$, and $c$ a tuple from $\C$, such that $a\in\dcl(Bc)$ (respectively, $a\in\acl(Bc)$).
\end{definition}

The type of any solution to $\delta(x)=1$, while being weakly $\C$-orthogonal, is at the same time $\C$-internal.

Examples of $\C$-orthogonality (or of non-almost-$\C$-internality, for that matter) are harder to verify as they require considering all possible extensions of parameters.
But examples do exist, even in order one: the type of any nonzero solution to $\delta(x)=\frac{x}{x+1}$ is orthogonal to the constants.
Such equations were studied by Kolchin, Rosenlicht, and Shelah; see Dave Marker's treatment of the above equation in~\cite[$\S$II.6]{mtof}.

Orthogonality to $\C$ implements the idea of having essentially no definable relation to the constants (even after passing to additional parameters), while almost $\C$-internality captures the opposite extreme of having a very significant definable relation with the constants (after possibly adding parameters).
There is, of course, a lot of room in-between.

Specialising to finite dimensional types, let us give the geometric meaning of $\C$-orthogonality and almost $\C$-internality for the generic types of $D$-varieties.
To do so we need a little more differential-algebraic geometry.
The following notions are very natural and I could have, probably should have, discussed them in the second lecture.

\begin{definition}
Suppose $(V_1,s_1)$ and $(V_2,s_2)$ are $D$-varieties over a differential field $(k,\delta)$.
Then a {\em $D$-rational map} $f:(V_1,s_1)\to(V_2,s_2)$ is a rational map $f:V_1\to V_2$ over $k$ such that $\tau f\circ s_1=s_2\circ f$.
That is, the diagram
$$\xymatrix{
\tau V_1\ar[r]^{\tau f}&\tau V_2\\
V_1\ar[u]^{s_1}\ar[r]^f&V_2\ar[u]_{s_2}
}$$
of rational maps, commutes.
\end{definition}

In the above definition I am using implicitly the fact that prolongation is a functor, so that $\tau f:\tau V_1\to\tau V_2$ makes sense.
When $f$ is a morphism it is quite clear how to define $\tau f$, namely: letting $\Gamma\subseteq V_1\times V_2$ be the subvariety given by the graph of $f$, verify that $\tau\Gamma\subseteq\tau(V_1\times V_2)=\tau V_1\times\tau V_2$ is in fact the graph of a morphism from $\tau V_1$ to $\tau V_2$, and then define $\tau f$ to be that morphism.
It is not too difficult to extend this construction to rational maps, but I leave the details to you.

The following is a useful way to check when a rational map is $D$-rational:

\begin{lemma}
\label{lem:drat}
Suppose $(V_i,s_i)$ are $D$-varieties over $(k,\delta)$ and $f:V_1\to V_2$ is a rational map.
Then the following are equivalent:
\begin{enumerate}
\item[(i)]
$f$ is $D$-rational,
\item[(ii)] $f\big((V_1,s_1)^\sharp\cap\dom(f)\big)\subseteq(V_2,s_2)^\sharp$, and
\item[(iii)]
$f(a)\in(V_2,s_2)^\sharp$ for some (equivalently any) generic $D$-point $a$ of $(V_1,s_1)$.
\end{enumerate}
\end{lemma}

\begin{proof}
I leave it to you to check that the following diagram
$$\xymatrix{
\tau V_1(\KK)\ar[r]^{\tau f}&\tau V_2(\KK)\\
V_1(\KK)\ar[u]^{\nabla}\ar[r]^f&V_2(\KK)\ar[u]_{\nabla}
}$$
always commutes.
From this, and the fact that $\nabla$ agrees with $s_i$ on $(V_i,s_i)^\sharp$, the implication (i)$\implies$(ii) follows easily.

Note that if $a$ is a generic $D$-point of $(V_1,s_1)$ over $k$ then it is Zariski-generic in $V_1$ over~$k$, and hence the rational map $f$ is defined at $a$.
Hence, that (ii) implies the ``for any" version of (iii) is clear.

Now, assume (iii) holds of some generic $D$-point $a$.
This implies that $\tau f\circ s_1$ and $s_2\circ f$ agree on $a$.
But $a$ is Zariski-generic in $\dom(f)$ over $k$, and such agreement is a Zariski closed condition over $k$.
It follows that they agree on all of $\dom(f)$, as required for~(i).
\end{proof}

OK, now back to orthogonality and internality for finite dimensional types.
Recall that the constants can be viewed as the set of $D$-points on the trivial vector field $(\AA^1,0)$.
So, from the differential-algebraic geometric point of view, a natural notion of ``interaction" between a $D$-variety $(V,s)$ and $(\AA^1,0)$ would be the existence of a nonconstant $D$-rational map $f:(V,s)\to(\AA^1,0)$.
By ``nonconstant" here we mean that $f$ is not a constant function on $V$, that as an element of the rational function field, $k(V)$, it is not in $k$.
However, the fact that $f$ is a $D$-rational map to $(\AA^1,0)$ does mean, exactly, that $f$ is a constant of the derivation $\delta_s$ on $k(V)$.
(You should check this.)
The constants of $(k(V),\delta_s)$ are often called {\em rational first integrals} for $(V,s)$.
In any case, here are the promised geometric characterisations:

\begin{theorem}
\label{thm:fdint}
Suppose $(V,s)$ is an absolutely irreducible $D$-variety over $(k,\delta)$  with generic type $p\in S(k)$.
Then
\begin{enumerate}
    \item[(a)] $p$ is nonorthogonal to the constants  if and only if there is a differential field extension $F\supseteq k$ and a nonconstant $D$-rational map $f:(V,s)\to(\AA^1,0)$ over $F$,
    \item[(b)] $p$ is almost internal to the constants if and only if there is a differential field extension $F\supseteq k$ and  a $D$-rational map $f:(V,s)\to(\AA^{\dim V},0)$ over $F$ which is dominant and generically finite-to-one.
    \end{enumerate}
\end{theorem}

\begin{proof}
Let us consider the right-to-left direction of part~(a).
Fix $a\in(V,s)^\sharp$ a generic $D$-point over $F$.
Then $a\models p$, and $a\ind_kF$ by Proposition~\ref{prop:indgentype}.
We have $f(a)\in (\AA^1,0)^\sharp=\C$.
Since $f$ is not constant on $V$, the Zariski locus of $a$ over $F\langle f(a)\rangle$ is a proper subset of $V$, and hence $a\nind_Ff(a)$, again by~\ref{prop:indgentype}.
This witnesses~$p\not\perp\C$.

Note that the same proof also gives the right-to-left direction of part~(b); indeed, in that case $f(a)\in\C^{\dim V}$ and $f$ being dominant and generically finite-to-one witnesses that $a\in\acl(Ff(a))$, which yeilds almost $\C$-internality.

For the converse of part~(a), let $L\supseteq k$ be a differential field extension and $c=(c_1,\dots,c_m)$ a tuple from $\C$ such that $a\ind_kL$ and $a\nind_Lc$.
That is, $L$ and $c$ witness that $p\not\perp\C$.
(Note that by taking definable closures we can always assume that the parameter extension witnessing nonorthogonality is a differential field extension.)
Moreover, suppose $m$ is minimal such.
So, if we set $F=L(c_1,\dots,c_{m-1})$, then $a\ind_kF$ and $a\nind_Fc_m$.
As $c_m$ is a constant, its dimension over any differential field is at most~$1$, and hence the fact that $c_m\nind_Fa$ implies that $c_m\in \acl(Fa)\setminus\acl(F)$.
Let $E$ be the finite orbit of $c_m$ under the action of the automorphisms of $(\KK,\delta)$ that fix $F$ and $a$ pointwise.
And let $e$ be a code for $E$.
(In a theory of fields, codes for finite sets always exist.)
So $e\in\dcl(Fa)\setminus\acl(F)$.
Moreover, as $E\subseteq\C$, we have that $e$ is itself a tuple of constants, say $e=(e_1,\dots,e_\ell)$.
Re-indexing, we may assume that $e_1\in\dcl(Fa)\setminus\acl(F)$.
Note that $\dcl(Fa)=F(a)$ and $a$ is Zariski-generic on $V$ over~$F$, so that $e_1=f(a)$ for some rational function $f$ on $V$ over~$F$.
That $f$ is nonconstant follows from the fact that $e_1\notin\acl(F)$.
Since $f(a)\in\C=(\AA^1,0)^\sharp$, Lemma~\ref{lem:drat} tells us that $f$ is a $D$-rational map, as desired.

The left-to-right direction of part~(b) is proved similarly.
Let $L\supseteq k$ and $c=(c_1,\dots,c_m)$ witness the almost $\C$-internality of $p$ so that $a\ind_kL$ and $a\in\acl(Lc)$.
Re-indexing we can find $0\leq \ell<m$ such that $(c_1,\dots,c_\ell)$ is an $\acl$-basis for $c$ over $La$.
Letting $F=L(c_1,\dots,c_\ell)$, we have $a\ind_kF$ and $a$ interalgebraic with $(c_{\ell+1},\dots,c_m)$ over~$F$.
Let $e$ be the code of the orbit of $(c_{\ell+1},\dots,c_m)$ over $Fa$ so that $e\in\dcl(Fa)$ and $a\in\acl(Fe)$.
So $e=f(a)$ where $f$ is a generically finite-to-one dominant rational map from $V$ to $W:=\zloc(e/F)$.
Exactly as before we see that $f:(V,s)\to(W,0)$ is $D$-rational.
As $\dim W=\dim(e/F)=\dim(a/F)=\dim(a/k)=\dim V$, we can compose $f$ with a finite-to-one co-ordinate projection dominantly onto $(\AA^{\dim V},0)$.
\end{proof}

\bigskip
\section{The dichotomy}

\noindent
In this final lecture it is my intention to articulate the Zilber dichotomy as it is manifest in $\dcf_0$, in terms of the birational geometry of $D$-varieties.
I will not prove anything here; neither the dichotomy itself, nor even that the geometric formulation I give is equivalent to the usual model-theoretic formulation.
The former (a proof of the dichotomy) is certainly beyond the scope of these lectures, while the latter I leave to the reader as it can be derived, with some work, using the various translations between model theory and geometry that we have discussed throughout these lectures.

The dichotomy will be a statement about $D$-varieties that are not ``covered by a family of proper infinite $D$-subvarieties".
More precisely:

\begin{definition}
\label{def:simple}
Suppose $(V,s)$ is an absolutely irreducible $D$-variety over $(k,\delta)$.
We say that $(V,s)$ is {\em simple}\footnote{My terminology is inspired by the bimeromorphic geometry of compact complex manifolds, where ``simplicity" is used to describe the analogous property.}
if the following holds:
for all irreducible $D$-varieties $(W,t)$ over $(k,\delta)$, and all proper irreducible $D$-subvarieties $Z$ of $(V\times W, s\times t)$ over~$k$, if $Z$ projects dominantly onto both~$V$ and~$W$ then $\dim Z=\dim W$.
\end{definition}

Let me explain a little how to think about this in terms of covering families of $D$-subvarieties.
Given an irreducible $D$-variety $(W,t)$ over $k$ and an irreducible $D$-subvarieties $Z$ of $(V\times W, s\times t)$ over~$k$ projecting dominantly onto~$W$, we view $Z$ as a {\em family of $D$-subvarieties of $(V,s)$ parametrised by $(W,t)$} in the following way: to each $b\in (W,t)^\sharp$ we can associate the fibre
$$Z_b:=\{x\in V:(x,b)\in Z\}.$$
It is a $D$-subvariety of $(V,s)$ over the differential field $k(b)$; this uses that $b$ is a $D$-point, $Z$ is a $D$-subvariety, and the co-ordinate projection is a $D$-morphism.
It is this family of $D$-subvarieties of $(V,s)$, namely $\{Z_b:b\in(W,t)^\sharp\}$, that we have in mind.
By a {\em generic member} of this family, or a {\em generic fibre}, we mean a fibre of the form $Z_b$ where $b$ is a generic $D$-point of $(W,t)$.
Note that the irreducibility of $Z$ ensures that this generic fibre $Z_b$ is irreducible over $k(b)$.
If $Z$ is a proper subvariety of $V\times W$ then the generic fibre $Z_b$ is a proper subvariety of $V$.
We say that $(V,s)$ is {\em generically covered} by the family if $Z$ projects dominantly onto $V$ as well.
The reason for this terminology is that it implies (indeed, precisely says) that if $a$ is a generic $D$-point of $(V,s)$ then $a\in Z_b$ for some (generic) $D$-point $b$ of $(W,t)$.
Note also that $\dim Z$ is the sum of $\dim W$ and the dimension of the generic fibre.
In particular, $\dim Z=\dim W$ if and only if the generic fibre is finite.
So Definition~\ref{def:simple} is saying that $(V,s)$ admits no generically covering family of $D$-subvarieties whose generic members are infinite and proper.

Model-theoretically, $(V,s)$ being simple is equivalent to the generic type of $(V,s)$ being {\em minimal}, that is, all its forking extensions are algebraic.
While we do not verify this here, let me say a few words that may be of use.
First of all, one has to convince oneself that the extra parameters needed to witness nonminimality can always be taken to be themselves finite dimensional.
So the generic $D$-points of the $(W,t)$ appearing in the definition of simple are the putative parameters for the forking extensions.
Secondly, if $\tp(a/k)$ is the generic type of $(V,s)$ and $\tp(a/kb)$ is a forking extension, with $\tp(b/k)$ the generic type of $(W,t)$, then the Zariski locus of $(a,b)$ over $k$ is a family of proper $D$-subvarieties of $(V,s)$ parametrised by $(W,t)$ that generically covers $(V,s)$.
Moreover, the converse holds as well; every such family gives rise to a forking extension.
Beyond that, one only has to verify that $\tp(a/kb)$ is algebraic if and only if the projection of $\zloc(a,b/k)$ onto $W$ is generically finite-to-one.

Every {\em $D$-curve}, by which I mean a $D$-variety $(V,s)$ where $\dim V=~1$, is simple.
This is because {\em every} proper subvariety of a curve is finite, let alone those that come from generically covering families.
But there are simple $D$-varieties of higher dimension.
Here is an example without proof:
Consider the surface $V\subseteq\AA^3$ defined by the equation $xz=1$ in co-ordinate variables $(x,y,z)$, and let $s$ be the section to the tangent bundle given by $s(x,y,z)=(x,y,z,y,yz,-yz^2)$.
Then $(V,s)$ is simple.
Indeed, this $D$-variety is the order~2 differential equation $xx''-x'=0$ in disguise, studied and shown to have minimal generic type by Poizat (see, for example, the treatment in~\cite[$\S$II.5]{mtof}).

In fact, simple $D$-varieties are ubiquitous in all dimensions.
Moreover, they are, in some real but subtle sense that I will not go into here, the building blocks for all $D$-varieties.
Suffice it to say that they are very much worth understanding.
The Zilber dichotomy will tell us that they either come from pure algebraic geometry or that they are geometrically very tame in a sense that we will now discuss.

\begin{definition}
Suppose $(V,s)$ is an absolutely irreducible $D$-variety over $(k,\delta)$.
A family $Z$ of $D$-subvarieties of $(V,s)$ parameterised by $(W,t)$ is said to be {\em rich} if the following conditions hold:
\begin{itemize}
\item[(i)]
for generic $b\in(W,t)^\sharp$, $Z_b$ is absolutely irreducible, and
\item[(ii)]
for generic $a\in(V,s)^\sharp$  there are infinitely many distinct generic fibres $Z_b$ that pass through $a$.
\end{itemize}
\end{definition}

Here by ``distinct generic fibres" I mean simply generic $D$-points $b_1,b_2$ of $(W,t)$ such that $Z_{b_1}\neq Z_{b_2}$ as subvarieties of $V$.

Note that in the case of a trivial vector field, every subvariety over constant parameters is a $D$-subvariety.
This is because the zero section takes every subvariety to its tangent bundle, which, when we are over constant parameters, agrees with the prolongation.
So trivial vector fields are a natural place to look for rich families of $D$-subvarieties, at least if the dimension is greater than~$1$.
For example, in $(\AA^2,0)$, we have the family of lines $y=mx+b$, which defines a  $D$-subvariety of $(\AA^2,0)\times(\AA^2,0)$ in the variables $(x,y,m,b)$, and hence a family of $D$-subvarieties of $(\AA^2,0)$ parameterises by $(\AA^2,0)$.
This family is rich because there are infinitely many lines through each point in the plane.

What the above example also illustrates, is that while simple $D$-varieties (more or less vaciously) admit no rich families of $D$-subvarieties, cartesian powers of a simple $D$-variety may: $(\AA^1,0)$ is simple but $(\AA^2,0)$ admits rich families of $D$-curves.
However, there are simple $D$-varieties all of whose cartesian powers admit no rich families: the Kolchin equation $x'=\frac{x}{x+1}$ and the Poizat equation $xx''=x'$, both of which we have already mentioned, are examples in dimension~$1$ and dimension~$2$ respectively.
The absence of rich families in all cartesian powers is a strong structural constraint on the differential-algebraic geometry of the $D$-variety.

The Zilber dichotomy states that $(\AA^1,0)$ is essentially the only simple $D$-variety that does not satisfy this strong structural constraint.

\begin{theorem}[Zilber dichotomy in $\dcf_0$, geometric formulation]
Suppose $(V,s)$ is a simple $D$-variety over $(k,\delta)$.
Then either
\begin{enumerate}
\item
there is a generically finite to one $D$-rational map $f:(V,s)\to(\AA^1,0)$ over some differential field extension of~$k$, or 
\item
whenever $n\geq 1$ and $X$ is an absolutely irreducible $D$-subvariety of $(V^n,s^n)$ over $k$ that projects dominantly onto $V$ in each co-ordinate, then $X$ admits no rich families.
\end{enumerate}
\end{theorem}

For the model theorist who is not used to seeing the dichotomy described this way, recall from Theorem~\ref{thm:fdint}(b) that case~(1) corresponds to the generic type of $(V,s)$ being almost $\C$-internal.
(Because of simplicity, in this case we actually get that $\dim V=1$.)
Case~(2) is equivalent to the generic type of $(V,s)$ being {\em $1$-based}, though I have not said anything about that.
In any case, this is the usual dichotomy in $\dcf_0$, but presented as a theorem about the birational geometry of $D$-varieties.

More is known about case~(2).
It splits into two more cases (thus forming the Zilber trichotomy); one coming from certain $D$-group structures on universal vectorial extensions of simple abelian varieties that do not descend to the constants (the {\em Manin kernel} case), and the other being when the cartesian powers of $(V,s)$ are truly devoid of any structure in that the absence of rich families is replaced by the absence of any infinite families at all (the {\em relationally trivial} case).
But that is a story for another course.

The Zilber dichotomy was first proved for $\dcf_0$  by Hrushovski and Sokolovic~\cite{HrSo}, relying on the theory of Zariski geometries developed by Hrushovski and Zilber~\cite{zg}.
Later, a new simpler proof was found by Pillay and Ziegler~\cite{PillayZiegler} using differential jet spaces, and having to do with the canonical base property, itself inspired by the model theory of compact complex manifolds.
But that too is a story for another course.

%\bibliography{differential}

\bigskip
\bigskip

\begin{thebibliography}{1}

\bibitem{Buium}
Alexandru Buium.
\newblock {\em Differential function fields and moduli of algebraic varieties}.
\newblock Springer-Verlag, 1986.

\bibitem{lamt}
Bradd Hart and Matthew Valeriote, editors.
\newblock {\em Lectures on algebraic model theory}, volume~15 of {\em Fields
  Institute Monographs}.
\newblock American Mathematical Society, Providence, RI, 2002.

\bibitem{HrSo}
Ehud Hrushovski and \v{Z}eljko Sokolovi\'{c}.
\newblock Strongly minimal sets in differentially closed fields.
\newblock {\em unpublished manuscript}, 1993.

\bibitem{zg}
Ehud Hrushovski and Boris Zilber.
\newblock Zariski geometries.
\newblock {\em Bull. Amer. Math. Soc. (N.S.)}, 28(2):315--323, 1993.

\bibitem{mtof}
David Marker, Margit Messmer, and Anand Pillay.
\newblock {\em Model theory of fields}, volume~5 of {\em Lecture Notes in
  Logic}.
\newblock Association for Symbolic Logic, La Jolla, CA; A K Peters, Ltd.,
  Wellesley, MA, second edition, 2006.

\bibitem{PillayZiegler}
Anand Pillay and Martin Ziegler.
\newblock Jet spaces of varieties over differential and difference fields.
\newblock {\em Selecta Math. (N. S.)}, 9 (4):579--599, 2003.

\end{thebibliography}
%\bibliographystyle{plain}

\end{document}